# THE MULTIPLIER ALGEBRA OF A NUCLEAR QUASIDIAGONAL $C^*$-ALGEBRA.

P. W. NG

ABSTRACT. The subject of quasidiagonality is of much interest in many places - among other things, in the classification program for simple unital separable nuclear $C^*$-algebras. In this note, we give two characterizations of nuclearity and quasidiagonality (for simple unital separable $C^*$-algebras).

Our first characterization is the nuclear analogue to Dadarlat's characterization of exact quasidiagonal $C^*$-algebras. Our second characterization is "dual" to the interesting (and important) popa property first studied by Popa.

## 1. INTRODUCTION

A unital $C^*$-algebra $\mathcal{A}$ is said to be *quasidiagonal* if there exists a faithful $*$-representation $\pi : \mathcal{A} \to \mathbb{B}(\mathcal{H})$ such that $\pi(\mathcal{A})$ is a set of quasidiagonal operators in $\mathbb{B}(\mathcal{H})$ - i.e., there exists an increasing sequence $\{p_n\}_{n=1}^\infty$ of finite rank projections in $\mathbb{B}(\mathcal{H})$ such that (a) $p_n$ converges to $1_{\mathbb{B}(\mathcal{H})}$ in the strong operator topology, and (b) for each $b \in \pi(\mathcal{A})$, $\|bp_n - p_n b\| \to 0$ as $n \to \infty$.

Note that in the definition in the previous paragraph, statements (a) and (b) are equivalent to saying that $\pi(\mathcal{A})$ is a set of (simultaneously) block diagonal operators in $\mathbb{B}(\mathcal{H})$, modulo the compacts. (The blocks are finite-dimensional.)

Notions around quasidiagonality play an interesting role in many places (see, for example, the survey paper [3] and the references therein), and have a long history with many contributions in different directions. In particular, quasidiagonality is an important notion in the classification program for simple unital separable nuclear $C^*$-algebras. Among other things, the following three questions are fundamental for classification theory (there are, of course, many other questions which follow from them):

**Question 1.1.** *Let $\mathcal{A}$ be a simple unital separable nuclear $C^*$-algebra.*

    i) *If $\mathcal{A}$ is stably finite then is $\mathcal{A}$ quasidiagonal?*

    ii) *If $\mathcal{A}$ is quasidiagonal then is $\mathcal{A}$ an inductive limit of type I $C^*$-algebras?*

    iii) *If $\mathcal{A}$ is quasidiagonal and satisfies some other condition - like $\mathcal{Z}$-stability or real rank zero - then can $\mathcal{A}$ be classified using $K$-theory invariants?? (In particular, in the real rank zero case, does $\mathcal{A}$ satisfy the tracial rank zero property of Huaxin Lin?)*

In connection with question i), we note that quasidiagonal $C^*$-algebras are stably finite (this follows from the definition by a short argument; see [3] proposition 3.19). In connection also with ii), all known classes of unital, simple, separable, nuclear, stably finite $C^*$-algebras that have been classified are quasidiagonal and are inductive limits of type I $C^*$-algebras (e.g. all simple unital $AH$-algebras with bounded dimension growth are certainly quasidiagonal). We also note that any nuclear quasidiagonal $C^*$-algebra is an inductive limit of nuclear residually finite dimensional $C^*$-algebras (see [3]). Finally, in connection with iii), an interesting recent discussion can be found in the appendix (section 13) of Brown's paper [5].

The notion of quasidiagonality has played a rather important role in the most recent work for classification in the stably finite case (with the other usual properties). As examples, we refer the reader to the papers [16], [6], [5], [14], [19] and [8].

Hence, a better understanding of nuclearity and quasidiagonality seems to be rather interesting. In this paper, we give (for simple unital separable $C^*$-algebras) two different characterizations of nuclearity and quasidiagonality - with hopes of applications to classification theory as well as norm topology amenability of the unitary group of a simple $C^*$-algebra. (The applications are in future papers and this paper basicly lays the groundwork).

Our first characterization of nuclearity and quasidiagonality (for simple unital separable $C^*$-algebras) is the nuclear analogue of Dadarlat's characterization of exact quasidiagonal $C^*$-algebras. In [7] (see also [3] Theorem 5.4), Dadarlat has given the following interesting characterization of separable exact quasidiagonal $C^*$-algebras:

**Theorem 1.2.** *Let $\mathcal{A}$ be a unital separable $C^*$-subalgebra. Then the following conditions are equivalent:*

i) $\mathcal{A}$ *is exact and quasidiagonal.*

ii) *Suppose that $\mathcal{A}$ is a unital $*$-subalgebra of $\mathbb{B}(\mathcal{H})$ such that $\mathcal{A} \cap \mathcal{K} = \{0\}$. Then for every finite subset $\mathcal{F} \subseteq \mathcal{A}$ and for every $\epsilon > 0$, there exists a finite-dimensional $C^*$-subalgebra $\mathcal{C} \subseteq \mathbb{B}(\mathcal{H})$ such that each element of $\mathcal{F}$ is within $\epsilon$ of an element of $\mathcal{C}$.*

Note that $\mathbb{B}(\mathcal{H})$ is actually the multiplier algebra of the algebra $\mathcal{K}$ of compact operators (i.e., $\mathbb{B}(\mathcal{H}) = \mathcal{M}(\mathcal{K})$). Using multiplier algebras, we give the simple nuclear analogue of Dadarlat's theorem:

**Theorem 1.3.** *Let $\mathcal{A}$ be a unital separable simple $C^*$-algebra. Then the following conditions are equivalent:*

i) $\mathcal{A}$ *is nuclear and quasidiagonal.*

ii) *Suppose that $\pi : \mathcal{A} \to \mathcal{M}(\mathcal{A} \otimes \mathcal{K})$ is a unital, purely large $*$-homomorphism. Then for every finite subset $\mathcal{F} \subseteq \pi(\mathcal{A})$ and for every $\epsilon > 0$, there exists a finite-dimensional $C^*$-subalgebra $\mathcal{C} \subseteq \mathcal{M}(\mathcal{A} \otimes \mathcal{K})$ such that each element of $\mathcal{F}$ is within $\epsilon$ of an element of $\mathcal{C}$.*

The *purely large property*, in theorem 1.3, is Elliott and Kucerovsky's concrete characterization of when an extension is nuclearly absorbing (see [11]; on the other hand, for the purposes of this paper, the reader will not be required to know the definition of "absorbing extension" and can ignore all statements containing this phrase!). From the absorbing extensions point of view, this property is the natural generalization (to arbitrary multiplier algebras) of the "essentialness" property (that $\mathcal{A} \cap \mathcal{K}(\mathcal{H}) = \{0\}$) in Dadarlat's theorem 1.2. This notion has also been quite important in classification theory. In section 2, we will describe this notion with some detail but the reader can skip many statements.

We also describe the *dual popa property* which is another characterization of nuclearity and quasidiagonality (for simple unital separable $C^*$-algebras), inspired by the work of Popa and Dadarlat (see section 4 for the precise definition of the dual popa property). The dual popa property is actually the opposite to the popa property studied by Popa (see [16] and [6]). Recall that Popa has given (for simple unital separable $C^*$-algebras with lots of projections (e.g. real rank zero)) a characterization of quasidiagonality which has been called by others the "popa property". This property inspired the tracial rank zero property of Lin (see [13]) which was an important step in the classification program.

One question in Popa's original paper ([16]) was whether or not simple unital real rank zero separable popa algebras with unique trace are always nuclear. Counterexamples to this (which are not hard to construct) were first discovered by Dadarlat and Lin. On the other hand, the dual popa property characterizes nuclearity and quasidiagonality (without needing to assume unique trace).

**Theorem 1.4.** *Let $\mathcal{A}$ be a unital separable simple $C^*$-algebra. Then the following conditions are equivalent:*

    i) *$\mathcal{A}$ is nuclear and quasidiagonal.*

    ii) *$\mathcal{A}$ has the dual popa property.*

The main technique used in this paper is from the theory of absorbing extensions, but the reader will not need to know anything about absorbing extensions! A good reference for extension theory is [1]. Good references for the theory of absorbing extensions (in the case of interest) are [11], [8], [12] and [1]. We note once more, though, that in this paper, the reader will not be required to know what an absorbing extension is! We also use some ideas from quasidiagonality. A good reference for quasidiagonality is Brown's survey paper [3] (and the references therein).

## 2. The purely large property

The *purely large property* is Elliott and Kucerovsky's concrete characterization of nuclearly absorbing extensions. (The reader will not need to know what an absorbing extension is and can freely ignore the statement of theorem 2.2 and other statements about absorbing extensions. But in the case of curiosity, please see [11] and the references therein.)

We here define the purely large property and state some results.

**Definition 2.1.** *Let* $\mathcal{A}$, $\mathcal{B}$ *and* $\mathcal{C}$ *be* $C^*$-*algebras, with* $B$ *stable. Suppose that we have an extension*

$$0 \to \mathcal{B} \to \mathcal{C} \to \mathcal{A} \to 0$$

*of* $\mathcal{B}$ *by* $\mathcal{A}$. *Then this extension is said to be* purely large *if for every positive element* $c \in C$ *which is not in* $\mathcal{B}$, *the* $C^*$-*algebra* $\overline{c\mathcal{B}c}$ *(the intersection of* $\mathcal{B}$ *with the hereditary subalgebra of* $\mathcal{C}$ *generated by c) is a stable full hereditary subalgebra of* $\mathcal{B}$.

Examples of purely large extensions are the extensions constructed by Kasparov and Lin (see, for example, [8]). These, as well as other purely large extensions, have played interesting roles in classification theory and other places (see, for example, [8], [11] etc...).

We also note that in the case where $\mathcal{B} = \mathcal{K}$ the $C^*$-algebra of compact operators (so the multiplier algebra is $\mathbb{B}(\mathcal{H})$), the "essentialness" condition $\mathcal{A} \cap \mathcal{K} = \{0\}$, in the statement of Dadarlat's result theorem 1.2, is equivalent to the purely large property.

Since we will be using these extensions, we now define (relevant special cases of) the Kasparov and Lin extensions.

**Definition 2.2.** *Let* $\mathcal{A}$ *be a simple unital separable* $C^*$-*algebra.*

(a) *Let* $\mathcal{H}$ *be a separable infinite dimensional Hilbert space, and let* $\pi :$ $\mathcal{A} \to \mathbb{B}(\mathcal{H})$ *be a unital* *-homomorphism. Then the extension, of* $\mathcal{A} \otimes \mathcal{K}$ by $\mathcal{A}$, given by the $*$-homomorphism $1_{\mathcal{A}} \otimes \pi : \mathcal{A} \to \mathcal{M}(\mathcal{A} \otimes \mathcal{K})$, is called a Kasparov extension.

(b) Let $id : \mathcal{A} \to \mathcal{A}$ be the identity map. The extension, of $\mathcal{A} \otimes \mathcal{K}$ by $\mathcal{A}$, given by the $*$-homomorphism $id \otimes 1_{\mathbb{B}(\mathcal{H})} : \mathcal{A} \to \mathcal{M}(\mathcal{A} \otimes \mathcal{K})$, is called the Lin extenstion.

**Proposition 2.1.** *Let $\mathcal{A}$ be a simple unital separable $C^*$-algebra. Then the Lin and Kasparov extensions, of $\mathcal{A} \otimes \mathcal{K}$ by $\mathcal{A}$ (as defined in the previous definition), are both purely large extensions.*

*Proof.* The proof involves applying the characterizations of stability developped by Hjelmborg and Rordam (see [17]). For details, see, among other things, [8]. □

The result about purely large extensions that we need (though the reader can skip), is the following result of Elliott and Kucerovsky [11]:

**Theorem 2.2.** *Let $\mathcal{A}$ and $\mathcal{B}$ be separable $C^*$-algebras, with $\mathcal{B}$ stable. Then an essential extension of $\mathcal{B}$ by $\mathcal{A}$ is nuclearly absorbing if and only if it is purely large.*

(Recall that an extension is said to be *essential* if the corresponding Busby map is one-to-one. Hence, the Lin and Kasparov extensions are both essential.)

Actually, for our purposes, we mainly need the following theorem from [11] (which follows from theorem 2.2):

**Theorem 2.3.** *Let $\mathcal{A}$ and $\mathcal{B}$ be separable $C^*$-algebras, with $\mathcal{A}$ unital and $\mathcal{B}$ stable. Suppose that either $\mathcal{A}$ or $\mathcal{B}$ is nuclear. Suppose that $\phi, \psi : \mathcal{A} \to \mathcal{M}(\mathcal{B})$ are two unital, essential and purely large $*$-homomorphisms. Then $\phi$ and $\psi$ are approximately unitarily equivalent in $\mathcal{M}(\mathcal{A} \otimes \mathcal{K})$.*

From this result and proposition 2.1, we have that (under the natural assumptions) the Lin and Kasparov extensions are approximately unitarily equivalent.

We end this section by pointing out that the theory of absorbing extensions play an important role in various places. For example, Kasparov used them to give a clean characterization of $KK$-theory (see [1]). They are also important in the stable existence and stable uniqueness theorems of classification theory (see [8]). Most recently, they have been used to classify certain interesting $C^*$-algebras associated with primitive aperiodic substitutional subshifts (dynamical systems). (See [9] and [10].)

## 3. The nuclear analogue of Dadarlat's Theorem

In this section, we prove the simple nuclear analogue of Dadarlat's characterization of exact quasidiagonal $C^*$-algebras.

**Theorem 3.1.** *Let $\mathcal{A}$ be a unital separable simple $C^*$-algebra. Then the following statements are equivalent:*

(1) *$\mathcal{A}$ is nuclear and quasidiagonal.*

(2) *If $\pi : \mathcal{A} \to \mathcal{M}(\mathcal{A} \otimes \mathcal{K})$ is a unital purely large $*$-homomorphism, then the image $\pi(\mathcal{A})$ can be locally approximated by finite-dimensional $*$-subalgebras of $\mathcal{M}(\mathcal{A} \otimes \mathcal{K})$.*

*Proof.* We first prove that (1) implies (2). Suppose that $\pi : \mathcal{A} \to \mathcal{M}(\mathcal{A} \otimes \mathcal{K})$ is a unital purely large $*$-homomorphism. Let $\epsilon > 0$ and a finite subset $\mathcal{F} \subseteq \mathcal{A}$ be given. Let $\gamma : \mathcal{A} \to \mathbb{B}(\mathcal{H})$ be a unital $*$-representation of $\mathcal{A}$. Let $\psi : \mathcal{A} \to \mathcal{M}(\mathcal{A} \otimes \mathcal{K})$ be the unital $*$-homomorphism given by $\psi(a) =_{df} 1_{\mathcal{A}} \otimes \gamma(a)$ for all $a \in \mathcal{A}$. Then $\psi$ is a Kasparov extension; and hence, by proposition 2.1 and theorem 2.3, $\pi$ and $\psi$ are approximately unitarily equivalent in $\mathcal{M}(\mathcal{A} \otimes \mathcal{K})$.

Since $\gamma$ is essential, it follows, from theorem 1.2 that $\gamma(\mathcal{A})$ can be locally approximated by unital finite-dimensional $*$-subalgebras of $\mathbb{B}(\mathcal{H})$. Hence, since $\pi$ and $\psi$ are approximately unitarily equivalent, $\pi(\mathcal{A})$ can be locally approximated by unital finite-dimensional $*$-subalgebras of $\mathcal{M}(\mathcal{A} \otimes \mathcal{K})$.

We now prove that (2) implies (1). We first prove nuclearity.

For each nonnegative integer $n$, let $\mathcal{A}_n$ be the $*$-subalgebra of $\mathcal{A}$ given by $\mathcal{A}_n =_{df} \{a \otimes 1_{\mathbb{B}(\mathcal{H})} : a \text{ is an element of } \mathbb{M}_{2^n}(\mathcal{A}). \}$. (Recall that $\mathcal{A} \otimes \mathcal{K} \cong \mathbb{M}_n(\mathcal{A}) \otimes \mathcal{K}$ for all $n$.) Conjugating by appropriate unitaries if necessary, we may assume that (a) $\{\mathcal{A}_n\}_{n=0}^{\infty}$ is an increasing sequence of $*$-subalgebras of $G$, and (b) $\bigcup_{n=0}^{\infty} \mathcal{A}_n$ is strictly dense in $\mathcal{A}$.

Now fix an integer $n$. Let $\epsilon > 0$ and a finite set of elements $\mathcal{G} \subseteq \mathcal{A}_n$ be given. The map $\mathcal{A} \to \mathcal{M}(\mathcal{A} \otimes \mathcal{K}) : a \mapsto a \otimes 1_{\mathbb{B}(\mathcal{H})}$ is a unital purely large $*$-homomorphism (this is the Lin extension and apply proposition 2.1). Hence, the image of this map can be locally approximated by finite-dimensional $*$-subalgebras of $\mathcal{M}(\mathcal{A} \otimes \mathcal{K})$. Hence, the image, of the map $\mathbb{M}_{2^n}(\mathcal{A}) \to \mathcal{M}(\mathcal{A} \otimes \mathcal{K}) : a \mapsto a \otimes 1_{\mathbb{B}(\mathcal{H})}$, can be locally approximated by finite-dimensional $*$-subalgebras of $\mathcal{M}(\mathcal{A} \otimes \mathcal{K})$. Hence, the elements of $\mathcal{G}$ are within $\epsilon$ of elements of a finite-dimensional $*$-subalgebra of $\mathcal{M}(\mathcal{A} \otimes \mathcal{K})$. But $\epsilon$ and $\mathcal{G}$ are arbitrary. Hence, since $\bigcup_{n=0}^{\infty} \mathcal{A}_n$ is strictly dense in $\mathcal{M}(\mathcal{A} \otimes \mathcal{K})$, $\mathcal{M}(\mathcal{A} \otimes \mathcal{K})$ can be locally strictly approximated by finite-dimensional $*$-subalgebras. Hence, since the strict topology is stronger than the $w*$-topology, $(\mathcal{A} \otimes \mathcal{K})^{**}$ is an injective von Neumann algebra. Hence, by [4], $\mathcal{A}$ is a nuclear $C^*$-algebra.

We now prove that $\mathcal{A}$ is quasidiagonal, by applying Voiculescu's abstract characterization of quasidiagonality. Once more, let $\epsilon > 0$ and a finite subset $\mathcal{F} \subseteq \mathcal{A}$ be given. Consider the unital $*$-homomorphism $\rho : \mathcal{A} \to \mathcal{M}(\mathcal{A} \otimes \mathcal{K}) : a \mapsto a \otimes 1_{\mathbb{B}(\mathcal{H})}$ (once more, the Lin extension). $\rho$ is a purely large extension, and hence, the image $\rho(\mathcal{A})$ can be locally approximated by finite-dimensional $*$-subalgebras of $\mathcal{M}(\mathcal{A} \otimes \mathcal{K})$. Choose a finite-dimensional $*$-subalgebra $F$ of $\mathcal{M}(\mathcal{A} \otimes \mathcal{K})$ such that $\rho(f)$ is "very norm-close" to an element of $F$ for all $f \in \mathcal{F}$. Since $F$ is a finite-dimensional $C^*$-algebra, it is an injective operator system. Hence, the identity map $F \to F$ extends to a unital completely positive map $\alpha : \mathcal{M}(\mathcal{A} \otimes \mathcal{K}) \to F$. Let $\beta : \mathcal{A} \to F$ be the map given by $\beta =_{df} \alpha \circ \rho$. $\beta$ is clearly a unital completely positive map. If $F$ was chosen to be a "good enough" local approximating algebra, then $\beta$ is also $\mathcal{F} - \epsilon$-multiplicative and $\mathcal{F} - \epsilon$-isometric. Since $\mathcal{F}$ and $\epsilon$ are arbitrary, it follows, by Voiculescu's abstract characterization of quasidiagonality ([3] theorem 4.2) that $\mathcal{A}$ is quasidiagonal.                                □

## 4. The dual popa property

For unital simple separable $C^*$-algebras, the *dual popa property* is a characterization of nuclearity and quasidiagonality. As indicated in the introduction, this property is the opposite of the popa property. Roughly speaking, a unital separable simple $C^*$-algebra has the popa property if given finitely many elements of $\mathcal{A}$, one can remove a possibly "big" piece (though never as big as the unit) and one will get a finite dimensional $C^*$-algebra (see [16] and [6]). The dual property, i.e., the dual popa property, says that given finitely many elements in $\mathcal{A}$, one can add a possibly "big" piece and one will get a finite dimensional $C^*$-algebra.

**Definition 4.1.** *Let $\mathcal{A}$ be a unital separable $C^*$-algebra. Then $\mathcal{A}$ is said to have the* dual popa property *if for every $\epsilon > 0$ and for every finite subset $\mathcal{F} \subseteq \mathcal{A}$, there is an integer $n$ and a completely positive $\mathcal{F} - \epsilon$-multiplicative map $\phi : \mathcal{A} \to \mathbb{M}_n(\mathcal{A})$ with $\phi(1_{\mathcal{A}})$ being a projection, there is a finite-dimensional $C^*$-subalgebra $\mathcal{D}$ of $\mathbb{M}_{n+1}(\mathcal{A})$, and there is a completely positive $\mathcal{F} - \epsilon$-multiplicative map $\Phi : \mathcal{A} \to \mathcal{D}$, such that $diag(f, \phi(f))$ is within $\epsilon$ of $\Phi(f)$ for every $f \in \mathcal{F}$.*

**Theorem 4.1.** *Let $\mathcal{A}$ be a unital separable simple $C^*$-algebra. Then the following conditions are equivalent:*

    (a) *$\mathcal{A}$ is nuclear and quasidiagonal.*

    (b) *$\mathcal{A}$ has the dual popa property.*

*Proof.* We first show that (b) implies (a). Let $\mathcal{F}$ be a finite subset of $\mathcal{A}$ and let $\epsilon > 0$ be given.

Since $\mathcal{A}$ has the dual popa property, choose maps $\Phi$ and $\phi$, and finite dimensional $*$-subalgebra $\mathcal{D}$ that satisfy the statements in definition 4.1. Now let $\psi : \mathbb{M}_{n+1}(\mathcal{A}) \to \mathcal{A}$ be the completely positive, completely contractive map gotten by taking the cutdown to the 1 by 1 entry. Since $\mathcal{D}$ is a subalgebra of $\mathbb{M}_{n+1}(\mathcal{A})$, we get the restricted map $\psi : \mathcal{D} \to \mathcal{A}$ (which we also denote by "$\psi$"). Moreover, $\psi \circ \Phi(f)$ is within $\epsilon$ of $f$, for all $f \in \mathcal{F}$.

Also, note that the map $\Phi$ is a completely positive, completely contractive, $\mathcal{F} - \epsilon$-multiplicative map that is also $\mathcal{F} - \epsilon$-isometric.

Since $\mathcal{F}$ and $\epsilon$ are arbitrary, it follows, from the definition of nuclearity, and from Voiculescu's abstract characterization of quasidiagonality (see [3] theorem 4.2), that $\mathcal{A}$ is nuclear and quasidiagonal.

The proof that (a) implies (b) follows from [7] proposition 2 and [2] proposition 6.1.6.

$\square$

## References

[1] B. Blackadar, *K-theory for operator algebras* Springer-Verlag, **New York** ( 1986) ;;

[2] B. Blackadar and E. Kirchberg, *Generalized inductive limits of finite-dimensional $C^*$-algebras* Math. Ann., **307** ( 1997) 343-380

[3] N. Brown, *On quasidiagonal $C^*$-algebras (in Operator algebras and applications 19-64)* Adv. Stud. Pure Math, **38** ( Math. Soc. Japan, Tokyo) 2004; A copy is available on the los Alamos server at arXiv:math.OA/0008181

[4] M. D. Choi and E. Effros, *Nuclear $C^*$-algebras and injectivity: the general case* Indiana Univ. Math. J., **261** ( 1977) 443-446

[5] N. Brown, *Invariant means and finite representation theory of $C^*$-algebras* A copy is available on the los Alamos server at http://xxx.lanl.gov/pdf/math.OA/0304009, () ;

[6] N. Brown, *Excision and a theorem of Popa* J. Operator Theory, **54** ( 2005) no. 1; 3-8; A copy is available at the los Alamos server at http://xxx.lanl.gov/pdf/math.OA/0305166;

[7] M. Dadarlat, *On the approximation of quasidiagonal $C^*$-algebras* J. Funct. Anal., **167** ( 1999) 69-78

[8] M. Dadarlat and S. Eilers, *On the classification of nuclear $C^*$-algebras* Proc. London Math. Society (3), **85** ( 2002) no. 1;168-210

[9] S. Eilers, G. Restorff and E. Ruiz, *lecture at the University of Muenster* Germany, ()

[10] S. Eilers, G. Restorff, and E. Ruiz, *lecture at the Fields Institute* Canada, ()

[11] G. A. Elliott and D. Kucerovsky, *An abstract Brown-Douglas-Fillmore absorption theorem* Pacific J. of Math., **3** ( 2001) 1-25

[12] H. Lin, *An introduction to the classification of amenable $C^*$-algebras* World Scientific Publishing Co., Inc., **River Edge, NJ** ( 2001) ;;

[13] H. Lin, *Classification of simple $C^*$-algebras with tracial topological rank zero* Duke Math. J., **125** ( 2004) no. 1; 91-119

[14] H. Lin, *Traces and simple $C^*$-algebras with tracial topological rank zero* J. Reine Angew. Math., **568** ( 2004) 99-137

[15] P. W. Ng, *Amenability of the sequence of unitary groups associated to a $C^*$-algebra* preprint, ()

[16] S. Popa, *On local finite-dimensional approximation of $C^*$-algebras* Pacific J. Math., **181** ( 1997) 141-158

[17] M. Rordam, *Stable $C^*$-algebras* Advanced studies in pure mathematics, **38** ( "Operator algebras and applications") edited by Hideki Kosaki; 2004; 195-264

[18] D. Voiculescu, *A note on quasidiagonal $C^*$-algebras and homotopy* Duke Math. J., **62** ( 1991) 267-271

[19] W. Winter, *On the classification of simple $Z$-stable $C^*$-algebras with real rank zero and finite decomposition rank* preprint, **A copy is available at the los Alamos archive at http://xxx.lanl.gov/pdf/math.OA/0502181** () ;

WestFaelische Wilhelms-Universitaet, Mathematisches Institut, Einstein-str. 62, 48149 Muenster, Germany, pwn@erdos.math.unb.ca

and current mailing address: The Fields Institute for research in the mathematical sciences, 222 College Street, Toronto, Ontario, M5T 3J1, Canada, pwn@erdos.math.unb.ca